\documentclass[11pt]{amsart}
\usepackage{a4,amsfonts,amssymb}
\usepackage{graphicx}
\usepackage{showlabels}

\theoremstyle{plain}
\newtheorem{thm}{Theorem}
\newtheorem{lem}{Lemma}

\theoremstyle{definition}
\newtheorem{df}{Definition}

\theoremstyle{remark}

\title[Delta shock wave and interactions  \dots]
{Delta shock wave and interactions in a simple model case}
\author{Marko Nedeljkov, Michael Oberguggenberger}

\begin{document}
\maketitle

\section{Introduction}
It has been observed by various authors \cite{hlf,k,kk,kor,mn1,tzz} that
the Riemann problem for certain equations from nonlinear elasticity and
gasdynamics cannot be solved for all combinations of piecewise
constant initial states with shock waves, rarefaction waves and contact
discontinuities only. For that reason, the notion of a delta shock wave
and a singular shock wave was introduced and employed by authors quoted
above, and it was shown that a large class of Riemann problems can be 
solved globally with these additional building blocks.
The aim of this paper is to study the interaction of one type of these new
solutions, the delta shock waves, with the classical types of solutions.

We continue the investigation of the model equation
\begin{align}
& u_{t}+(u^{2}/2)_{x}=0 \label{e1}\\
& v_{t}+((u-1)v)_{x}=0 \label{e2}
\end{align}
initiated in \cite{hlf}.
This system is derived from a simplified model of magneto-hydrodynamics. 
In \cite{hlf}, the authors 
found a solution for every Riemann problem with the initial data 
$(u_{0},v_{0})$ on the left- and $(u_{1},v_{1})$ on the right-hand side 
from zero in the following way.

The eigenvalues of the above system are $\lambda_{1}(u,v)=u-1$,
$\lambda_{2}(u,v)=u$, and the right-hand side eigenvectors are 
$r_{1}(u,v)=(0,1)^{T}$, $r_{2}(u,v)=(1,v)^{T}$. The first characteristic
field is linearly degenerate and the second is genuinely nonlinear. Thus, 
there are three types of solution.

\noindent
(i)
When $u_{1}>u_{0}$ the solution is a contact discontinuity followed by a 
rarefaction wave,
\begin{align*}
u(x,t) & = \begin{cases} u_{0}, & x\leq u_{0}t \\
{x\over t}, & u_{0}t<x<u_{1}t \\
u_{1}, & x\geq u_{1}t \end{cases} \\
v(x,t) & = \begin{cases} v_{0}, & x\leq (u_{0}-1)t \\
v_{1}\exp(u_{0}-u_{1}), & (u_{0}-1)t<x<u_{0}t \\
v_{1}\exp({x\over t}-u_{1}), & u_{0}t\leq x \leq u_{1}t \\
v_{1}, & x>u_{1}t. \end{cases}
\end{align*}
 
\noindent (ii) If $u_{1}<u_{0}<u_{1}+2$, the solution is given in the 
form of contact discontinuity followed by a shock wave,
\begin{align*}
u(x,t) & = \begin{cases} u_{0}, & x\leq ct \\
u_{1}, & x> ct \end{cases} \\
v(x,t) & = \begin{cases} v_{0}, & x\leq (u_{0}-1)t \\
v_{\ast}, & (u_{0}-1)t<x<ct \\
v_{1}, & x \geq ct, \end{cases}
\end{align*}  
where $\displaystyle v_{\ast}=v_{1}{2-u_{0}-u_{1} \over 2+u_{1}-u_{0}}$.  

\noindent (iii) If $u_{0}\geq u_{1}+2$ the solution is given in the form 
of delta shock wave, 
\begin{align*}
u(x,t) & = \begin{cases} u_{0}, & x\leq ct \\
u_{1}, & x> ct \end{cases} \\
v(x,t) & = \left\{ \begin{aligned} v_{0}, & \;\; x\leq ct \\
v_{1}, & \;\; x> ct \end{aligned}\right\}+\alpha_{0}(t)D^{-}+\alpha_{1}(t)D^{+},
\end{align*}
where $D^{-}$ and $D^{+}$ are the left- and right-hand side 
delta functions with the support on the line $x=ct$ (see below),
$c=(u_{0}+u_{1})/2$,
\begin{equation*}
\alpha_{0}(t)={st(c-(u_{1}-1)) \over u_{0}-u_{1}},\;\; 
\alpha_{1}(t)={st(c-(u_{0}-1)) \over u_{0}-u_{1}},
\end{equation*}
$\alpha(t):=\alpha_{0}(t)+\alpha_{1}(t)$ is called the strength of
the delta shock wave, and
\begin{equation*}
s:=c(v_{1}-v_{0})-((u_{1}-1)v_{1}-(u_{0}-1)v_{0})
\end{equation*}
is called the Rankine-Hugoniot deficit (see \cite{kk}).
\medskip \medskip

Our aim is to investigate various possible interactions of a solution 
in one of these forms with a delta shock wave. There are five possibilities
for this to happen.
\begin{itemize}
\item[ Case 1.] delta shock wave interact with an another one 
\item[ Case 2.] delta shock wave interact
with a contact discontinuity followed by a shock wave from the left-hand side 
\item[ Case 3.] delta shock wave interact with a contact discontinuity
followed by a shock wave from the right-hand side 
\item[ Case 4.] delta shock wave interact with a contact discontinuity
followed by a rarefaction wave from the left-hand side 
\item[ Case 5.] delta shock wave interact with a contact discontinuity
followed by a rarefaction wave from the right-hand side 
\end{itemize}
We shall always assume that the shock wave or the rarefaction wave starts
from $(0,0)$ and the delta shock wave from another point to left or right from 
zero. The initial data are determined by triplets 
$(u_{0},u_{1},u_{2})$ and $(v_{0},v_{1},v_{2})$. 

We shall now briefly describe what we mean by a solution in the form
of a delta shock wave.

Suppose $\overline{R_{+}^{2}}$ is divided into finitely disjoint
open sets $\Omega_{i}\neq \emptyset$, $i=1,...,n$ with piecewise 
smooth boundary curves $\Gamma_{i}$, $i=1,...,m$, that is 
$\Omega_{i}\cap \Omega_{j}=\emptyset$, $\bigcup_{i=1}^{n}
\overline{\Omega}_{i}=\overline{R_{+}^{2}}$ where $\overline{\Omega}_{i}$
denotes the closure of $\Omega_{i}$. Let 
${\mathcal C}(\overline{\Omega}_{i})$ be the space of bounded
and continuous real-valued functions on $\overline{\Omega}_{i}$,
equipped with the $L^{\infty}$-norm. Let ${\mathcal M}(\overline{\Omega}_{i})$,
be the space of measures on $\overline{\Omega}_{i}$.

We consider the spaces
\begin{equation*}
{\mathcal C}_{\Gamma}=\prod_{i=1}^{n}{\mathcal C}(\overline{\Omega}_{i}),
\quad 
{\mathcal M}_{\Gamma}=\prod_{i=1}^{n}{\mathcal M}(\overline{\Omega}_{i}).
\end{equation*}
The product of an element $G=(G_{1},...,G_{n})\in {\mathcal C}_{\Gamma}$ and
$D=(D_{1},...,D_{n})\in {\mathcal M}_{\Gamma}$ is defined as an element
$D\cdot G=(D_{1}G_{1},...,D_{n}G_{n})\in {\mathcal M}_{\Gamma}$,
where each component is defined as the usual product of a continuous
function and a measure.

Every measure on $\overline{\Omega}_{i}$ can be viewed as a measure on 
$\overline{{\mathbb R}_{+}^{2}}$ with support in $\overline{\Omega}_{i}$.
This way we obtain a mapping
\begin{equation*}
\begin{split}
& m:{\mathcal M}_{\Gamma} \rightarrow {\mathcal M}
(\overline{{\mathbb R}_{+}^{2}}) \\
& m(D)=D_{1}+D_{2}+...+D_{n}.
\end{split}
\end{equation*}
A typical example is obtained when $\overline{{\mathbb R}_{+}^{2}}$
is divided into two regions $\Omega_{1}$, $\Omega_{2}$ by a piecewise
smooth curve $x=\gamma(t)$. The delta function $\delta(x-\gamma(t))
\in {\mathcal M}(\overline{{\mathbb R}_{+}^{2}})$ 
along the line $x=\gamma(t)$ can 
be split in a non unique way into a left-hand side  
$D^{-}\in {\mathcal M}(\overline{\Omega}_{1})$ and the right-hand component
$D^{+}\in {\mathcal M}(\overline{\Omega}_{2})$ such that
\begin{equation*}
\begin{split}
\delta(x-\gamma(t)) & = \alpha_{0}(t)D^{-}+\alpha_{1}(t)D^{+} \\
& = m(\alpha_{0}(t)D^{-}+\alpha_{1}(t)D^{+})
\end{split}
\end{equation*}
with $\alpha_{0}(t)+\alpha_{1}(t)=1$.
The solution concept which allows to incorporate such two sided delta
functions as well as shock waves is modeled along the lines of the 
classical weak solution concept and proceeds as follows:

\noindent
Step 1: Perform all nonlinear operations of functions in the space
${\mathcal C}_{\Gamma}$.

\noindent
Step 2: Perform multiplications with measures in the space 
${\mathcal M}_{\Gamma}$.

\noindent
Step 3: Map the space ${\mathcal M}_{\Gamma}$ into 
${\mathcal M}(\overline{{\mathbb R}_{+}^{2}})$ by means of the map $m$
and embed it into the space of distributions.

\noindent
Step 4: Perform the differentiation in the sense of distributions and require
that the equation is satisfied in this sense.
\medskip

Note that in the case of absence of a measure part (Step 2), this is 
the precisely the concept of a weak solution to equations in divergence form.
\medskip

Following the reasoning in \cite{kk}, delta shocks are required to satisfy
the condition of overcompressibility, meaning that all characteristic
curves run into the delta shock curve from both sides. It may happen 
that at a certain point on a delta shock curve, overcompressibility 
is lost. In this case we replace the delta shock by a new type of solution
which we call a delta contact discontinuity. This new concept is introduced
in Lemma 1 and Definition 1 below.

At interaction points our solutions are computed by continuation as 
continuous functions of time with values in the space of distributions
(as solutions to a new initial value problem at the time of interaction).

The result of our investigation  is that the interaction of a delta shock wave
with any of three types of solutions to the Riemann problem, (i)--(iii)
above, can be described by means of delta shocks and delta contact
discontinuities. This is summarized in the following theorem.

\begin{thm}\label{Theorem}
The initial value problem for system (\ref{e1}, \ref{e2}) with three
constant states, one of which produces a delta shock, has a global weak
solution consisting of a combination of rarefaction waves, shock waves,
contact discontinuities, delta shock waves and delta contact 
discontinuities.
\end{thm}

The remainder of the paper is devoted to proving 
this result by going throug all possible cases of interaction.

\section{Interactions with shock waves}

\noindent
{\bf Case 1.}
Here $u_{0}\geq u_{1}+2$, $u_{1}\geq u_{2}+2$. The speeds of the 
delta shock waves are $c_{1}=(u_{0}+u_{1})/2$ and $c_{1}=(u_{1}+u_{2})/2$.
At the interaction point $(x_{0},t_{0})$, the new initial data are
\begin{equation*}
\begin{split}
& u|_{t=t_{0}}=\begin{cases} u_{0}, & x<x_{0} \\ u_{2}, & x>x_{0} \end{cases}\\
& v|_{t=t_{0}}=\left\{ \begin{aligned} v_{0}, & \;\; x<x_{0} \\ v_{2}, 
& \;\; x>x_{0} \end{aligned} \right\}
+\gamma \delta_{(x_{0},t_{0})},
\end{split}
\end{equation*}
where $\gamma$ denotes a sum of the strengths of incoming delta shock waves.

Let $u=G$, $v=H+(\alpha_{0}(t)D^{-}
+\alpha_{1}(t)D^{+})$, where $G$ and $H$ are step functions
\begin{equation*}
G=\begin{cases}u_{0}, & x-x_{0}<(t-t_{0})c \\
u_{2},& x-x_{0}>(t-t_{0})c \end{cases}
\quad  H=\begin{cases}v_{0}, & x-x_{0}<(t-t_{0})c \\
v_{2},& x-x_{0}>(t-t_{0})c, \end{cases}
\end{equation*}
and $D=\alpha_{0}(t)D^{-}+\alpha_{1}(t)D^{+}$ is a split delta
function supported by the line $x=x_{0}+(t-t_{0})c$.

From (\ref{e1}) it follows
\begin{equation*}
-c[G]+{1\over 2}[G^{2}]=0,
\end{equation*}
i.e.\
$c=(u_{0}+u_{2})/2$. Since $c_{1}>c>c_{2}$, the wave will be the 
overcompressive one, because of $u_{0}-1>c>u_{2}$. Equation (\ref{e2}) gives
\begin{equation*}
\begin{split}
& -c[H]\delta+[(G-1)H]\delta+(\alpha_{0}'(t)+\alpha_{1}'(t))\delta \\
& -c(\alpha_{0}(t)+\alpha_{1}(t))\delta'+
((u_{0}-1)\alpha_{0}(t)+(u_{2}-1)\alpha_{1}(t))\delta'=0
\end{split}
\end{equation*}
This equation gives (with $\alpha(t)=\alpha_{0}(t)+\alpha_{1}(t)$))
the following ODE
\begin{equation*}
\alpha'(t)=c(v_{2}-v_{0})-((u_{2}-1)v_{2}
-(u_{0}-1)v_{0})=:s\in {\mathbb R},\;\; \alpha(t_{0})=\gamma.
\end{equation*}
The unique solution is given by $\alpha(t)=s(t-t_{0})+\gamma$. Substitution
of $\alpha$ in the equation gives 
\begin{equation*}
\begin{split}
& \alpha_{0}(t)+\alpha_{1}(t)=\alpha(t)=s(t-t_{0})+\gamma \\
& (u_{0}-1-cs(t-t_{0})-\gamma)\alpha_{0}(t)
+(u_{2}-1-cs(t-t_{0})-\gamma)\alpha_{1}(t)=0 
\end{split}
\end{equation*}
which has a unique solution 
$\alpha_{0}(t)$, $\alpha_{1}(t)$, since $u_{0}\neq u_{2}$.

Thus, the result of the first type of interaction is a single 
delta shock wave.
\medskip 

\noindent
{\bf Case 2.}
Suppose that the delta shock wave is given by
\begin{equation*}
\begin{split}
& u(x,t)=G(x+a^{2}-c_{1}t), \\ 
& v(x,t)=H(x+a^{2}-c_{1}t)
+\beta_{0}(t)\delta^{-}(x+a^{2}-c_{1}t)
+\beta_{1}(t)\delta^{+}(x+a^{2}-c_{1}t),\\
&G=\begin{cases} u_{0}, & x+a^{2}<c_{1}t \\
u_{1}, & x+a^{2}>c_{1}t \end{cases},\quad 
H=\begin{cases} v_{0}, & x+a^{2}<c_{1}t \\
v_{1}, & x+a^{2}>c_{1}t, \end{cases}
\end{split}
\end{equation*}
and that a contact discontinuity coupled with a shock wave is given by
\begin{equation*}
\begin{split}
&u=\begin{cases} u_{1}, & x<c_{2}t \\
u_{1}, & x>c_{2}t \end{cases},\quad 
v=\begin{cases} v_{1}, & x<(u_{1}-1)t \\
v_{\ast}, & (u_{1}-1)t\leq x<c_{2}t \\ 
v_{2}, & x\geq c_{2}t, \end{cases}
\end{split}
\end{equation*}
where $c_{2}=(u_{1}+u_{2})/2<c_{1}$ (since $u_{1}>u_{2}$, $u_{0}\geq u_{1}+2$),
and $v_{\ast}=v_{2}(2+u_{1}-u_{2})/(2+u_{2}-u_{1})$.

Let us denote by $(t_{0},x_{0})$ the point where delta shock wave meets
the contact discontinuity, i.e.\
this point is the intersection of the lines $x+a^{2}=c_{1}t$ and 
$x=(u_{1}-1)t$.

In the area bounded by the lines $x=(u_{1}-1)t$ and $x=u_{1}t$, the value
of $u$ is the constant $u_{1}$. This implies that the delta shock wave runs 
through it with the same speed $c_{1}=(u_{0}+u_{1})/2$ as before. Only 
the values of $\beta_{0}(t)$ and $\beta_{1}(t)$ are changed into, say, 
$\tilde{\beta}_{0}(t)$ and $\tilde{\beta}_{1}(t)$ due to the existing 
difference in $v_{0}$ and $v_{\ast}$. The new strength of the delta shock
wave is now $s_{1}(t-t_{0})+\gamma_{0}$, where 
\begin{equation*}
s_{1}:=c_{1}(v_{\ast}-v_{0})-(u_{0}-1)(v_{\ast}-v_{0}),
\end{equation*}
and $\gamma_{0}$ is the strength of the previous delta shock wave in the 
point $(t_{0},x_{0})$. Obviously, the new delta shock wave is 
an overcompressive wave, since $u_{0}-1\geq c_{1} \geq u_{1}$.

Let us denote by $(t_{1},x_{1})$ the point where the new delta shock 
meets the existing shock wave i.e.\
the point $(t_{1},x_{1})$ is the intersection of the lines 
$x+c_{1}^{2}=c_{1}t$ and $x=c_{2}t$. Let $\gamma_{1}$ be the strength 
of the delta shock wave at this point. Therefore, we obtain the new initial
data 
\begin{equation*}
u=\begin{cases} u_{0}, & x<x_{1} \\ u_{2}, & x>x_{1} \end{cases} \quad 
v=\left \{ \begin{aligned} v_{0}, & \;\; x<x_{1} \\ v_{2}, & \;\; x>x_{1} 
\end{aligned}\right\} 
+\gamma_{1}\delta_{(t_{1},x_{1})}.
\end{equation*}  
A solution for the new initial data problem will be a delta shock wave
with the speed $c=(u_{0}+u_{2})/2<c_{1}$, and $c$ is obtained directly from 
(\ref{e1}) in the usual way. Again this speed ensures that the obtained 
wave is an overcompressive one, since $u_{0}-1\geq c\geq u_{2}$. 

Substituting $u$ and $v$ into (\ref{e2}) gives the strength
\begin{equation*}
\tilde{\alpha}_{0}(t)+\tilde{\alpha}_{1}(t)=s_{1}(t-t_{1})+\gamma_{1}.
\end{equation*}
Using this equation and 
\begin{equation*}
(u_{0}-1-s(t-t_{1})-\gamma_{1})\tilde{\alpha}_{0}(t)
+(u_{2}-1-s(t-t_{1})-\gamma_{1})\tilde{\alpha}_{1}(t)=0,
\end{equation*}
where $s:=c(v_{2}-v_{0})-((u_{2}-1)v_{2}-(u_{0}-1)v_{0})$, one can 
find unique $\tilde{\alpha}_{0}(t)$ and $\tilde{\alpha}_{1}(t)$, 
and this proves the above statement.
\medskip 

\noindent
{\bf Case 3.} Now, $u_{0}>u_{1}>u_{2}+2$, and the speed of the shock wave 
$c_{1}=(u_{0}+u_{1})/2$ is greater than the speed of the delta shock wave
$c_{2}=(u_{1}+u_{2})/2$. Let $(t_{0},x_{0})$ be the interaction point 
of these two waves, and let $\gamma_{0}$ be the strength of the delta shock
wave at this point. Initial data are now
\begin{equation*}
u|_{t=0}=\begin{cases} u_{0}, & x<x_{0} \\
u_{2}, & x>x_{0} \end{cases} \quad 
v|_{t=0}=\left\{ \begin{aligned} v_{\ast}, & \;\; x<x_{0} \\
v_{2}, & \;\; x>x_{0} \end{aligned}\right\}+\gamma_{0}\delta_{(t_{0},x_{0})},
\end{equation*}
where the value of $v_{\ast}$ is defined as before.

Similarly to the previous case, the result of the interaction is a single
overcompressive delta shock wave with the speed $c=(u_{0}+u_{2})/2$ 
(since $u_{0}-1\geq c\geq u_{2}$). As before, $c$ is obtained from (\ref{e1})
and  from (\ref{e2}) one can find $\tilde{\alpha}_{0}(t)$ 
and $\tilde{\alpha}_{1}(t)$ in the same way as above.

All the way through the contact discontinuity, the delta shock wave has the 
same speed, only $\tilde{\alpha}_{0}(t)$ and $\tilde{\alpha}_{1}(t)$ 
are changing.

\section{Interactions with rarefaction waves}

One can easily see that interaction of a rarefaction and delta shock wave 
are much more complicated. Now we shall deal with this problem.

As one could see before, the new initial data include a delta function 
as a part. If the right-hand side of $u$ is greater
or equal to the left-hand one plus 2, the new initial value problem can be 
solved in a simple way as above
and the result is a single overcompressive delta
shock wave. But, when this is not a case, the types of  
admissible solution known so far are not enough to obtain a solution.
The definition of a new type of admissible solution,
called delta contact discontinuity, is given below.  
Its existence is justified by two facts. 
First, a contact discontinuity emerges in the case when one of the 
characteristic fields is linearly degenerate. Second, if a linear equation has
a delta function as initial data, it propagates along the 
characteristic lines. These two facts inspired the following lemma
and the definition of this new type of elementary waves. 
 
\begin{lem} \label{l1}
Let the initial data for system (\ref{e1}-\ref{e2}) given by 
\begin{equation*}
u|_{t=0}=\begin{cases} u_{0}, &  x<0 \\ u_{1}, & x>0 \end{cases},\;\;
v|_{t=0}=\left\{ \begin{aligned} v_{0}, & \;\; x<0 \\ v_{1}, & \;\; x>0 
\end{aligned} \right\} +\gamma \delta_{(0,0)},
\end{equation*}
where $u_{0}>u_{1}$, but $u_{0}<u_{1}+2$. Then, the function
\begin{equation*}
u=\begin{cases} u_{0}, & x<ct \\ u_{1}, & x>ct \end{cases}\;\;
v=\left\{ \begin{aligned}  v_{0}, & \;\; x<(u_{0}-1)t \\ v_{\ast}, 
& \;\; (u_{0}-1)t<x<ct \\
v_{1}, & x>ct \end{aligned} \right\} + \gamma\delta_{x=(u_{0}-1)t},
\end{equation*}
where $c=(u_{0}+u_{1})/2$ weakly solves the Riemann problem for 
(\ref{e1},\ref{e2}). 
\end{lem}
\begin{proof}
For every $\varphi\in C_{0}^{\infty}$, $\mathop{\rm supp}\varphi
\cap \{(x,t):\;\; x=(u_{0}-1)t, \;\; t>0\}=\emptyset$, it holds that
\begin{equation*}
\begin{split}
& \langle u_{t},\varphi\rangle 
+{1\over 2} \langle (u^{2})_{x},\varphi\rangle=0 \\
& \langle v_{t},\varphi\rangle 
+ \langle ((u-1)v)_{x},\varphi\rangle=0. 
\end{split}
\end{equation*}
Our aim is to show that this still holds true
when it is allowed that 
$\mathop{\rm supp}\varphi$ intersects the supports of $D^{-}$ and $D^{+}$, i.e.\
the line $x=(u_{0}-1)t$. 
Let us note that the condition $u_{0}<u_{1}+2$ means that 
$(u_{0}+u_{1})/2>u_{0}-1$ so the line $x=(u_{0}-1)t$ is on the left-hand side
of the shock line $x=(u_{0}+u_{1})t/2$.

Equation (\ref{e1}) does not contain $v$, so it is still 
satisfied. From (\ref{e2}) we have that 
\begin{equation*}
\begin{split}
v_{t}+((u_{0}-1)v)_{x}&=-(u_{0}-1)(v_{\ast}-v_{0})\delta-\gamma\delta' \\
& +(u_{0}-1)(v_{\ast}-v_{0})\delta+(u_{0}-1)\gamma \delta' =0 
\end{split}
\end{equation*}
near the line $x=(u_{0}-1)t$.
\end{proof}

Usefulness of this lemma will be clear after the interaction of a
delta shock and rarefaction wave is treated. Then one could roughly see 
how a solution looks like, since the rare faction wave could be approximated
with a large number of small amplitude non-physical shock waves
(see \cite{Bre}, for example).

Another possible use could be in a sort of
a wave front tracking algorithm, where systems in question posses a solution
containing a delta function.

\begin{df}
Consider a region $R$ where $u$ is continuous function and a curve 
$\Gamma$ in $R$ of slope $\lambda_{1}(u,v)$. A distribution 
$(u,v)\in {\mathcal C}(R)\times {\mathcal D}'(R)$ is
a delta contact discontinuity, if $v$ is a sum of a locally integrable
function on $R$ and a delta function on $\Gamma$ which weakly solves
(\ref{e1},\ref{e2}) on $R$.
\end{df}

Let us note that the overcompressiveness condition obviously need not  
hold in this case. But, as we already have mentioned before the lemma, 
a linearly degenerate field resembles a linear equation where
this type of a solutions exists.
\medskip 

Also, we shell try to show admissibility of the delta contact discontinuity
using the entropy and entropy-flux functions for system (\ref{e1}-\ref{e2}).
Entropy and appropriate
entropy-flux functions are given by
\begin{equation*}
\begin{split}
\eta(u,v)= & f(u)+g(v e^{-u})e^{u} \\
q(u,v)= & e^{u}g(v e^{-u})u-e^{u}g(v e^{-u})+\int uf'(u) du \\
= & (u-1)\eta(u,v)+\tilde{f}(u).
\end{split}
\end{equation*}
Substituting the functions $u$ and $v$ in a neighbourhood of the
delta contact discontinuity support $x=(u_{0}-1)t$ by piecewise
constant functions
$$ u\equiv u_{0},\; 
v=\begin{cases} v_{0}, & x<(u_{0}-1)t-\varepsilon \\
v_{1 \varepsilon}, & (u_{0}-1)t-\varepsilon <x< (u_{0}-1)t \\
v_{2 \varepsilon}, & (u_{0}-1)t <x< (u_{0}-1)t +\varepsilon \\
v_{\ast}, & (u_{0}-1)t+\varepsilon <x
\end{cases}, $$
where $\gamma=\lim_{\varepsilon \rightarrow 0}\varepsilon
(v_{1 \varepsilon}+v_{2 \varepsilon})$.
one gets
$$ \eta(u,v)_{t}+q(u,v)_{x}\approx 0.$$
That is, convex entropy condition is satisfied for each entropy 
function pair. 
\medskip \medskip

Now, we are returning to the last two cases which covers the rest of
possible delta shock wave interactions.

\noindent
{\bf Case 4.}  Suppose that a delta shock wave starts from the 
point $-a^{2}$, $a>0$, with the speed $c_{1}=(u_{0}+u_{1})/2$ and
meets a contact discontinuity followed by
the rarefaction wave 
centered at zero. Denote by $(\tilde{t}_{0},\tilde{x}_{0})$ the
meeting point, i.e.\
it is the intersection of the lines $x+a^{2}=c_{1}t$ and 
$x=(u_{1}-1)t$. As we have already seen, the delta shock wave goes
through the contact discontinuity without speed change 
(but its strength is changed) and meets
the rarefaction wave at some point $(x_{0},t_{0})$,
\begin{equation*}
t_{0}={2a^{2} \over u_{0}-u_{1}},\;\; x_{0}={2u_{1}a \over u_{0}-u_{1}}.
\end{equation*}
Let $\gamma_{0}$ be the strength of the delta shock wave at this point.
In order to see what could happen, let us approximate
the rarefaction wave with a set of non-physical shock waves, supported
by the lines $x=(u_{1}+n\eta)t$, $\eta<<1$, $n\in {\mathbb N}$.
(see Fig.~1.)
 
\begin{center}
\includegraphics*[scale=0.65]{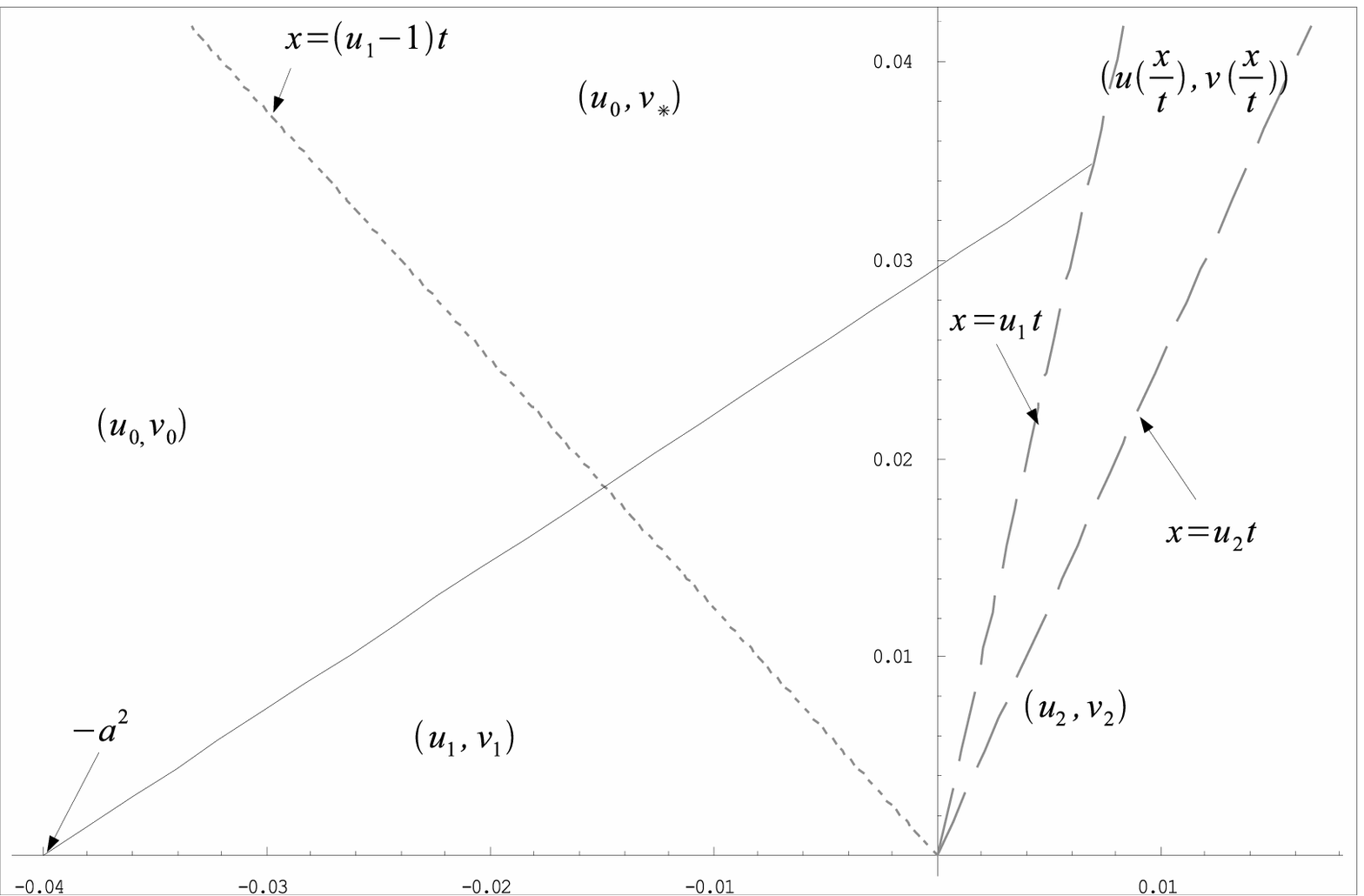}

{\small Fig.~1.}
\end{center}

At least in the beginning, until $(u_{1}+\eta n)+2\leq u_{0}$, 
the result of successive interactions of the delta shock wave with
the non-physical shock waves are delta shock waves with increasing
speeds, with values $(u_{0},v_{0})$ on the left-hand side and 
the values on the right-hand side are the values of the rarefaction wave.
This guide us to look for a curve $\Gamma_{0}:=(c(t),t)$, such that 
a delta function lives on it, $c(t_{0})=x_{0}$. The value of 
$u$ on the left-hand side of $\Gamma$ is $u_{0}$, and $c(t)/t=x/t$ on the 
right-hand side. Inserting the above data for such a curve into 
(\ref{e1}), one gets the following ordinary differential equation
\begin{equation}\label{Eq}
-c'(t)\Big( {c(t)\over t}-u_{0}\Big)+
{1\over 2} \Big( \Big( {c(t)\over t}\Big)^{2}-u_{0}^{2}\Big)=0,
\;\; c(t_{0})=x_{0},
\end{equation}
which has the unique solution
\begin{equation}\label{e3}
c(t)=u_{0}t-a\sqrt{2(u_{0}-u_{1})t},\;\; t\geq t_{0}.
\end{equation}
Denote by $v(t)$ the value of $v|_{\Gamma_{0}}$ in the rarefaction wave,
$v(t)=v_{2}\exp(c(t)/t-u_{2})$. Substituting expected delta shock wave
given by $u=G$, $v=H+\alpha_{0}(t)D_{\Gamma_{0}}^{-}+
\alpha_{1}(t)D_{\Gamma_{0}}^{+}$, where $G$ and $H$ are the step 
functions with discontinuity line $\Gamma_{0}$, gives
\begin{equation*}
\begin{split}
& -c'(t)(v(t)-v_{\ast})\delta+c(t)(\alpha_{0}(t)+\alpha_{1}(t))'\delta \\
& -c'(t)(\alpha_{0}(t)+\alpha_{1}(t))\delta'+
\left( \left({c(t) \over t}\right) v(t)-(u_{0}-1)v_{\ast}\right)\delta \\
& + \left( (u_{0}-1)\alpha_{0}(t)+\left( {c(t)\over t}-1\right)
\alpha_{1}(t)\right)\delta'=0.
\end{split}
\end{equation*}

Since the following ordinary differential equation
\begin{equation*}
\alpha'(t)={c'(t)(v(t)-v_{\ast})-(u_{0}-1)v_{\ast}+\big( {c(t)\over t}-1\big)
v(t) \over c(t)},\;\; \alpha(t_{1})=\gamma_{1}
\end{equation*}
has a unique solution (obtained in a simple manner by an integration),
the strength of the delta shock wave, $\alpha(t)$, is determined. 

Equating the coefficient of $\delta'$ with zero, we can compute the
two summands $\alpha_{0}$ and $\alpha_{1}$ of $\alpha$.
Since 
\begin{equation*}
c'(t)=u_{0}-{a\sqrt{2(u_{0}-u_{1})}\over 2\sqrt{t}}
>{c(t) \over t}=u_{0}-{a\sqrt{2(u_{0}-u_{1})}\over \sqrt{t}},
\end{equation*}
the obtained delta shock wave satisties rhe right-hand overcompressibility
condition. Overcompressibility
condition for the left-hand side is 
\begin{equation}\label{equ4}
u_{0}-1\geq c'(t).
\end{equation}
Now, we have the following two cases.

\noindent
(i) If $u_{2}\leq u_{0}-2$, relation (\ref{equ4}) is satisfied trough all the 
rarefaction wave and the resulting solution is a single delta
shock wave with the speed $c=(u_{0}+u_{2})/2$ starting from the
point $(\tilde{x},\tilde{t})$ which is the intersection 
of the curve $\Gamma_{0}$ and the line $x=u_{2}t$.

\begin{center}
\includegraphics*[scale=0.65]{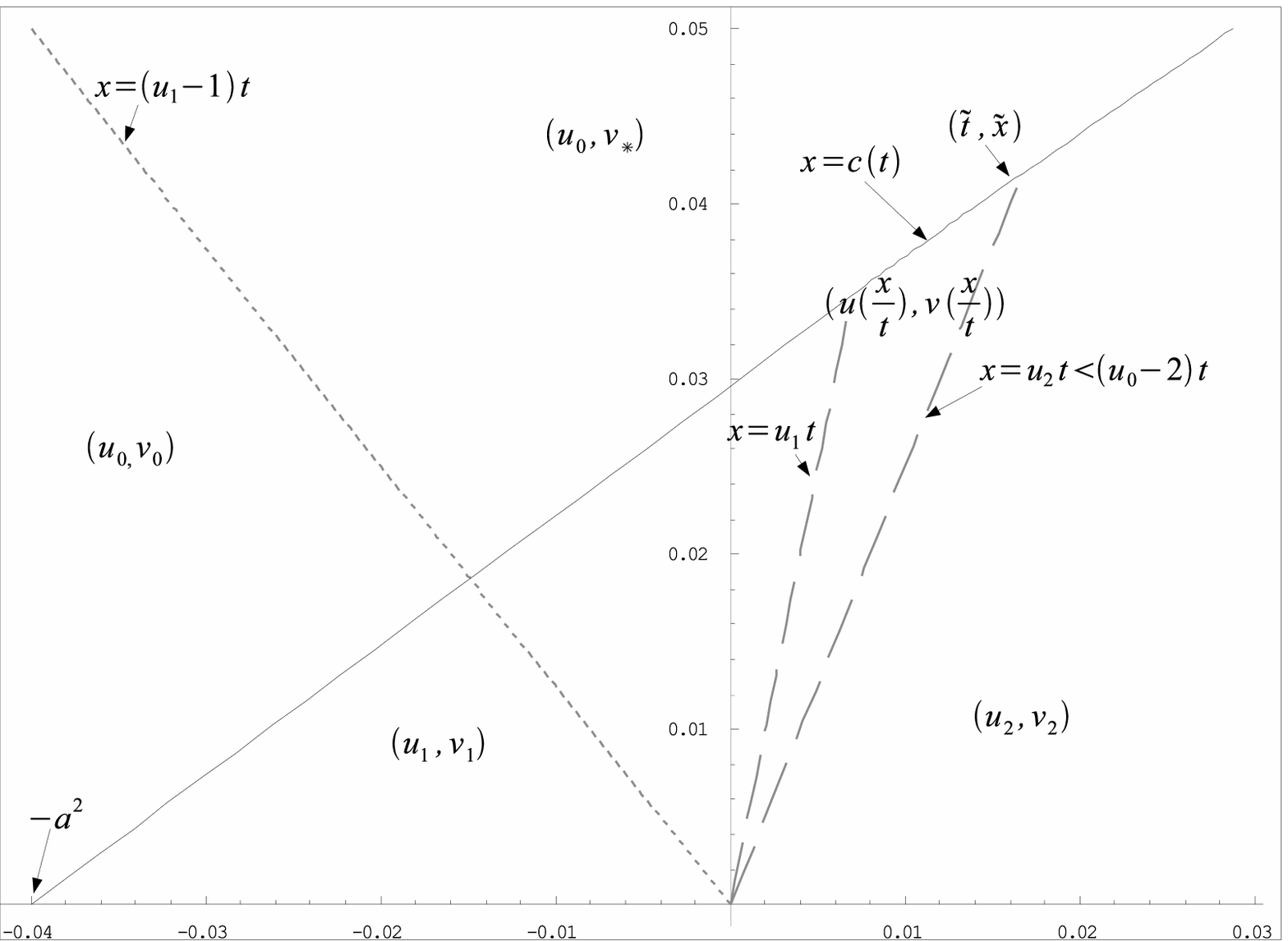}

{\small Fig.~2.}
\end{center}

After the time $\tilde{t}$, the solution in this case is given by 
\begin{equation*}
u|_{t>\tilde{t}}=\begin{cases} u_{0}, &   
x<\tilde{x}t \\
u_{2}, & x>\tilde{x}t \end{cases} \quad
v|_{t>\tilde{t}}=\left \{ \begin{aligned} v_{\ast}, & \;\; 
x<\tilde{x}t \\
v_{2},&  \;\; x>\tilde{x}t \end{aligned}\right\} + \gamma_{(\tilde{x}.\tilde{t})}
\delta_{(\tilde{x},\tilde{t})} 
\end{equation*}

\noindent
(ii) Suppose that $u_{2}>u_{0}-2$. Then the delta shock wave supported
by $\Gamma_{0}$ is an overcompressive wave only until some 
point $(x_{s},t_{s})$ lying inside the rarefaction wave.
(See Fig.~3.)

\begin{center}
\includegraphics*[scale=0.65]{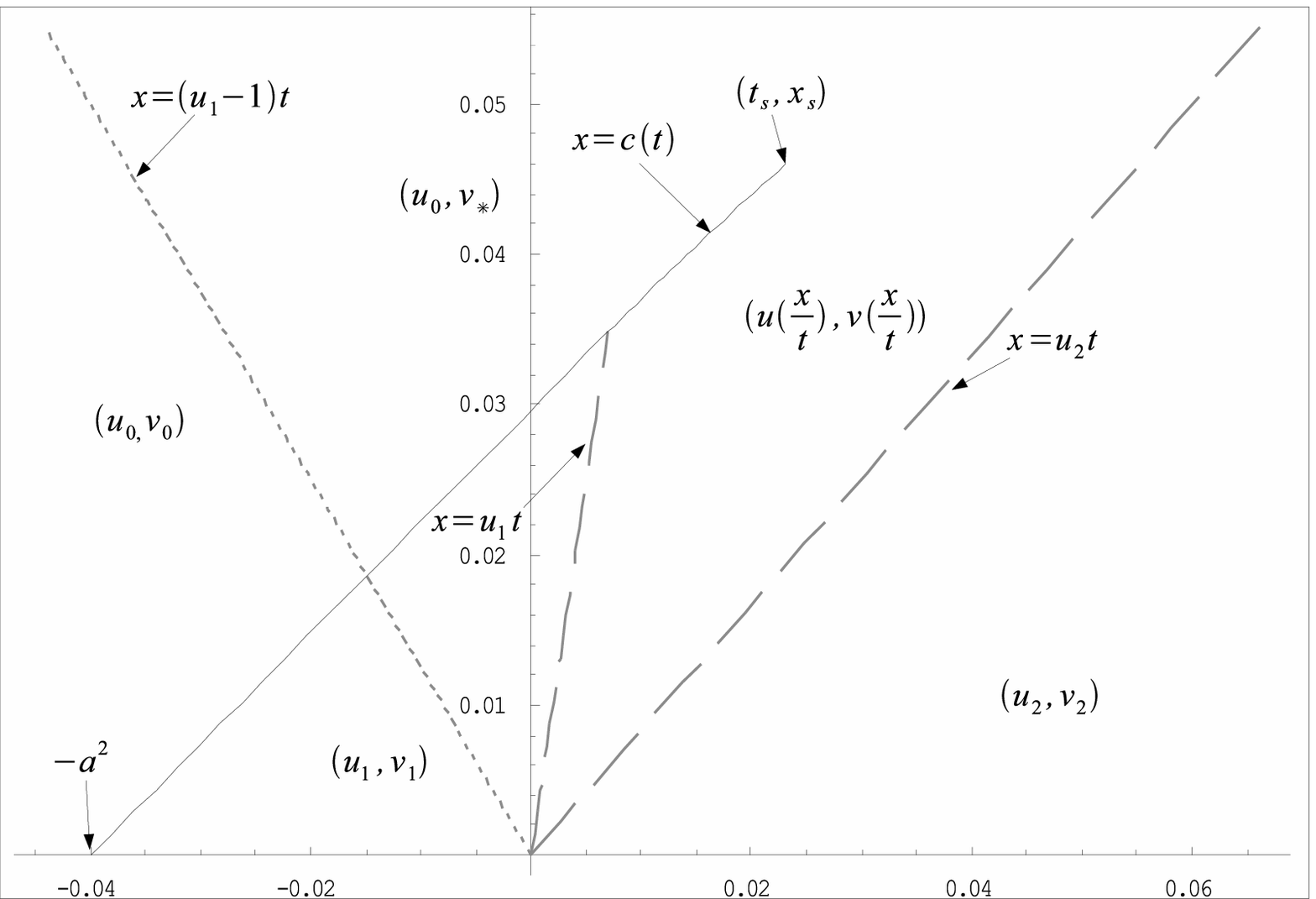}

{\small Fig.~3.}
\end{center}

So, the admissible solution cannot be prolonged along the same curve 
$\Gamma_{0}$. Assuming that the rarefaction wave is approximated
by a set of small non-physical shock waves, the present problem is 
described in Lemma \ref{l1}: the right-hand side
equals $u_{0}-2+\eta$, $0<\eta\ll 1$, while the left-hand
one equals $u_{0}$.

In this lemma, the problem is solved by using the new type of 
a solution -- delta contact discontinuity. This is exactly what
we shall try. That is, suppose that the solution consists of the delta
function supported by a line $\Gamma_{1}:\; (x-x_{s})=(u_{0}-1)(t-t_{s})$
going through an area where $u$ has a constant value $u_{0}$, and 
a shock wave supported by a curve $\Gamma_{2}:\;\; x=c_{2}(t)$, 
where $c_{2}(t_{s})=x_{s}$, with the left-hand side values 
$u_{0}$ of the function $u$ and the right-hand side ones 
$c_{2}(t)/t$ (a part of the rarefaction wave). All that means
that $c_{2}(t)$ should satisfy the same equation (\ref{Eq}) as $c(t)$
with the initial data $c_{2}(t_{s})=x_{s}=c(t_{s})$, i.e.\
the new shock wave is supported by the continuation of the  
curve $\Gamma_{0}$.

\begin{center}
\includegraphics*[scale=0.65]{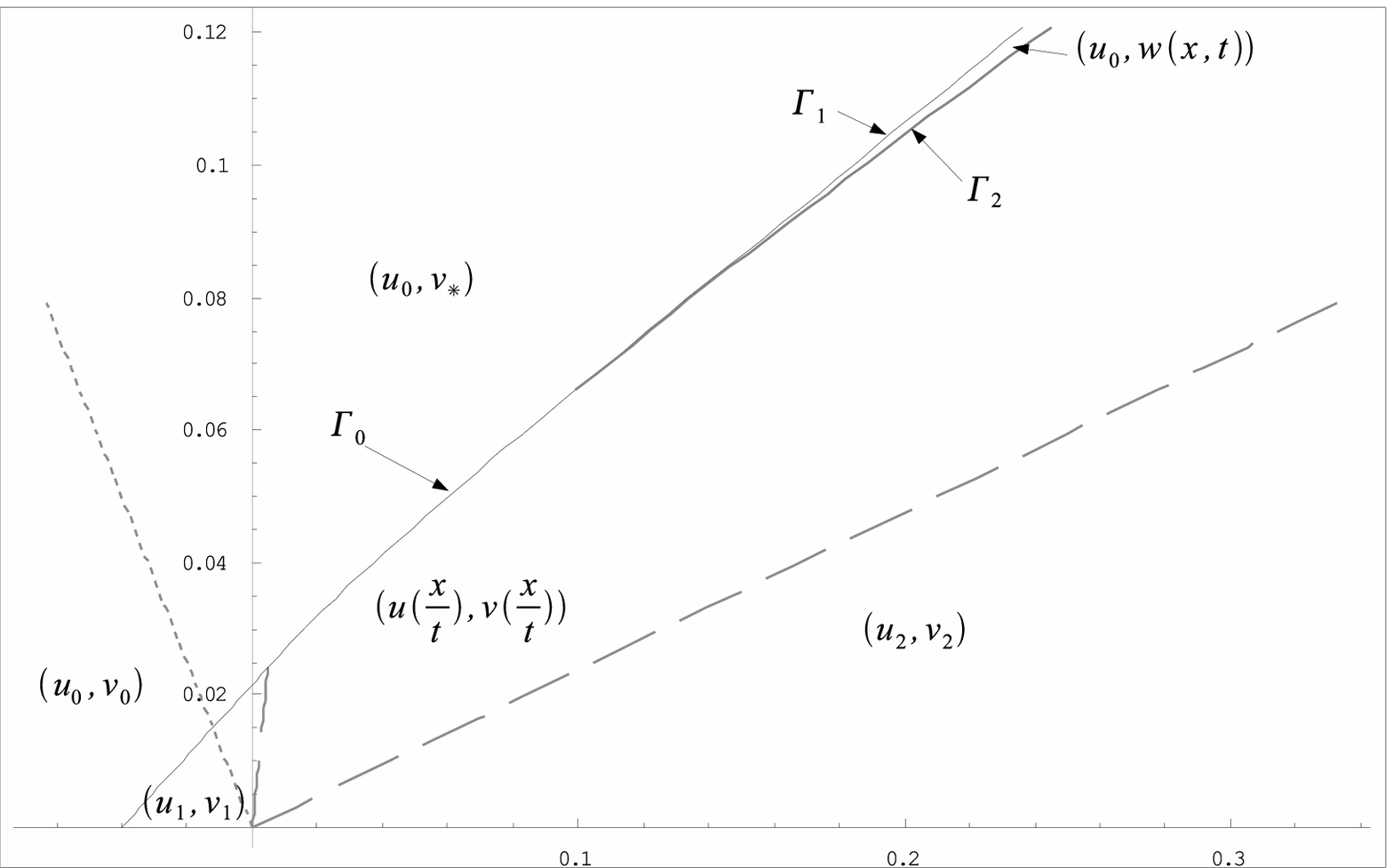}

{\small Fig.~4.}
\end{center}

Since $u_{0}>c_{1}'(t)$ and $u_{0}-1<c_{1}'(t)$ while 
$c_{1}'(t)>c(t)/t$, the obtained shock wave, supported by the
curve $\Gamma_{2}$ is admissible. 

The Rankine-Hugoniot conditions for $\Gamma_{2}$ after the time
$t=t_{s}$ imply
\begin{equation} \label{e5}
-c_{2}'(t)(v(t)-w_{\ast}(t))+\left( {c_{2}(t) \over t}-1\right)
v(t)-(u_{0}-1)w_{\ast}(t)=0,
\end{equation}
where $w_{\ast}$ denotes the left-hand side value of $v$ along the curve
$\Gamma_{2}$. Equation (\ref{e5}) simply determines
\begin{equation*}
\begin{split}
w_{\ast}(t)&= {\sqrt{t}+B \over \sqrt{t}-B}v(t) \\
& = {\sqrt{t}+B \over \sqrt{t}-B}v_{2}\exp(u_{0}t-2B\sqrt{t}-u_{1}),
\end{split}
\end{equation*}
where $B:=a\sqrt{2(u_{0}-u_{1})}/2=\sqrt{t_{s}}$.

The value of $v$ between $\Gamma_{1}$ and $\Gamma_{2}$, 
denoted by $w(x,t)$ has to satisfy
the  equation
\begin{equation}\label{e6}
w_{t}+(u_{0}-1)w_{x}=0, \;\; w_{\Gamma_{1}}=w_{\ast}(t).
\end{equation}
The solution to (\ref{e6}) is of the form $w(x,t)=V(y)$, $y=x-(u_{0}-1)t$.
More precisely, using the initial data one gets
\begin{equation} \label{e7}
V(y)=v_{2}\left(1+{2B \over \sqrt{B^{2}+y}}\right)
\exp((u_{0}B+u_{0}\sqrt{B^{2}+y}-2B)(B+\sqrt{B^{2}+y})-u_{1}).
\end{equation}
The curve $\Gamma_{1}$ is given by
\begin{equation} \label{e8}
x-(u_{0}-1)t=x_{s}-(u_{0}-1)t_{s}=(u_{s}-u_{0}-1)t_{s}=-t_{s}.
\end{equation}
Substitution of (\ref{e8}) into (\ref{e7}) yields $V|_{\Gamma_{2}}=\infty$
and $V(y)\in {\mathbb R}$ for $(t,x)$ lying between $\Gamma_{1}$ and 
$\Gamma_{2}$, since $B^{2}=t_{s}$. But $w\in L_{\mathop{\rm loc}}^{1}
\subset {\mathcal D}'$. 

In order to verify that it is a solution we note that
\begin{equation*}
v(t,x)=\left\{ \begin{aligned} v_{\ast}, & \;\; x<(u_{0}-1)t \\ 
w(x,t),& \;\; x>(u_{0}-1)t, \;\; x<c_{2}(t) 
\end{aligned} \right\}
+\gamma_{s}\delta_{\Gamma_{2}},
\end{equation*}
when $\gamma_{s}\delta_{\Gamma_{2}}$ is the delta function with the 
strength $\gamma_{s}$ obtaining from the initial data at $(t_{s},x_{s})$.
Since $v(x,t)$ is constant along the lines parallel to $x=(u_{0}-1)t$
in a region where $u\equiv u_{0}$ it is clear that it is a solution of
(\ref{e2}).

In order to see what is going on after the interaction
of the delta shock and rarefaction wave, one has to consider three 
different possibilities.
\medskip 

\noindent
(a) $u_{0}\leq u_{2}$.

Then the delta contact discontinuity and shock wave supported by 
$\Gamma_{1}$ lies inside the rarefaction wave since 
$\Gamma_{1} \cap \{(t,x):\;\; x=u_{2}t\}=\emptyset$ and
$\Gamma_{2} \cap \{(t,x):\;\; x=u_{2}t\}=\emptyset$ 
(actually, $c_{2}(t)$ has the line $x=u_{0}t$ as an asymptote, as
$t \rightarrow \infty$ (see \ref{e3})).
\medskip

\noindent
(b) $u_{0}>u_{2}\geq u_{0}-1$.

Then the delta contact discontinuity stays inside the rarefaction wave 
and the shock wave supported by $\Gamma_{2}$ intersects the line
$x=u_{2}t$ at some point $(\tilde{t},\tilde{x})$.

Now, at the point $(\tilde{t},\tilde{x})$ we have the new Cauchy
problem for (\ref{e1},\ref{e2}) with the initial data $(u_{0},w(x,t))$,
$(u_{2},v_{2})$. Since $u_{0}>u_{2}$, but $u_{0}<u_{2}+2$, the 
solution is given by 
\begin{equation*}
\begin{split}
& u=\begin{cases} u_{0}, & x-\tilde{x}<(u_{0}+u_{2})(t-\tilde{t})/2 \\
u_{2}, & x-\tilde{x}>(u_{0}+u_{2})(t-\tilde{t})/2 \end{cases} \\
& v=\begin{cases} w(t,x), & x-\tilde{x}<(u_{0}-1)(t-\tilde{t}) \\
\tilde{v}_{\ast}, 
& (u_{0}-1)(t-\tilde{t})<x-\tilde{x}<(u_{0}+u_{2})(t-\tilde{t})/2 \\
v_{2}, & x-\tilde{x}>(u_{0}+u_{2})(t-\tilde{t})/2, \end{cases} 
\end{split}
\end{equation*}
where  
$\tilde{v}_{\ast}=v_{2}(2+u_{0}-u_{2})/(2+u_{2}-u_{0})$. Let us remark that
the function $w$ equals a constant value along lines with the slope
$(u_{0}-1)$. Denote by $\Gamma_{4}$ the shock line 
$x-\tilde{x}=(u_{0}+u_{2})(t-\tilde{t})/2$. This line is a tangent 
to the curve $\Gamma_{2}$ at $(\tilde{x},\tilde{t})$.

Let $\Gamma_{3}$ be the line of slope $u_{0}-1$ starting at 
$(\tilde{t},\tilde{x})$.
Since $\tilde{t}$ is a solution to 
\begin{equation*}
u_{0}t-2B\sqrt{t}=u_{2}t,
\end{equation*}
we have 
\begin{equation*}
\tilde{t}={4B^{2} \over (u_{0}-u_{2})^{2}}, \;\; 
\tilde{x}={4u_{2}B^{2} \over (u_{0}-u_{2})^{2}}
\end{equation*}
and 
\begin{equation*}
w|_{\Gamma_{3}}={\sqrt{\tilde{t}}+B \over \sqrt{\tilde{t}}-B}v(\tilde{t})
={2+u_{0}-u_{2} \over 2+u_{2}-u_{0}}v_{2}=\tilde{v}_{\ast},
\end{equation*}
after the substitution of $\tilde{t}$ and the ending value
of the rarefaction wave, $v(\tilde{t})=v_{2}$. So, the function
$w(t,x)$ is continuously prolonged by $\tilde{v}_{\ast}$ into the area 
between the lines $x-\tilde{x}=(u_{0}-1)(t-\tilde{t})$ 
(the contact discontinuity line) and 
$x-\tilde{x}=(u_{0}-u_{2})(t-\tilde{t})/2$ 
(the shock curve).

The slope of $\Gamma_{3}$ is the same as the 
one of $\Gamma_{1}$. That is, there are no interactions,
and this case is finished.

\noindent 
(c) $u_{2}<u_{0}-1$. In this case both of $\Gamma_{1}$ and $\Gamma_{2}$
intersects the line $x=c_{2}t$. But as $\Gamma_{2}$ reaches this line 
at the time $\tilde{t}=4B^{2}/(u_{0}-u_{2})^{2}$, before the time 
when $\Gamma_{1}$ would intersect it, analysis is the same as in the 
case (b) (see Fig.~5) below.

\begin{center}
\includegraphics*[scale=0.65]{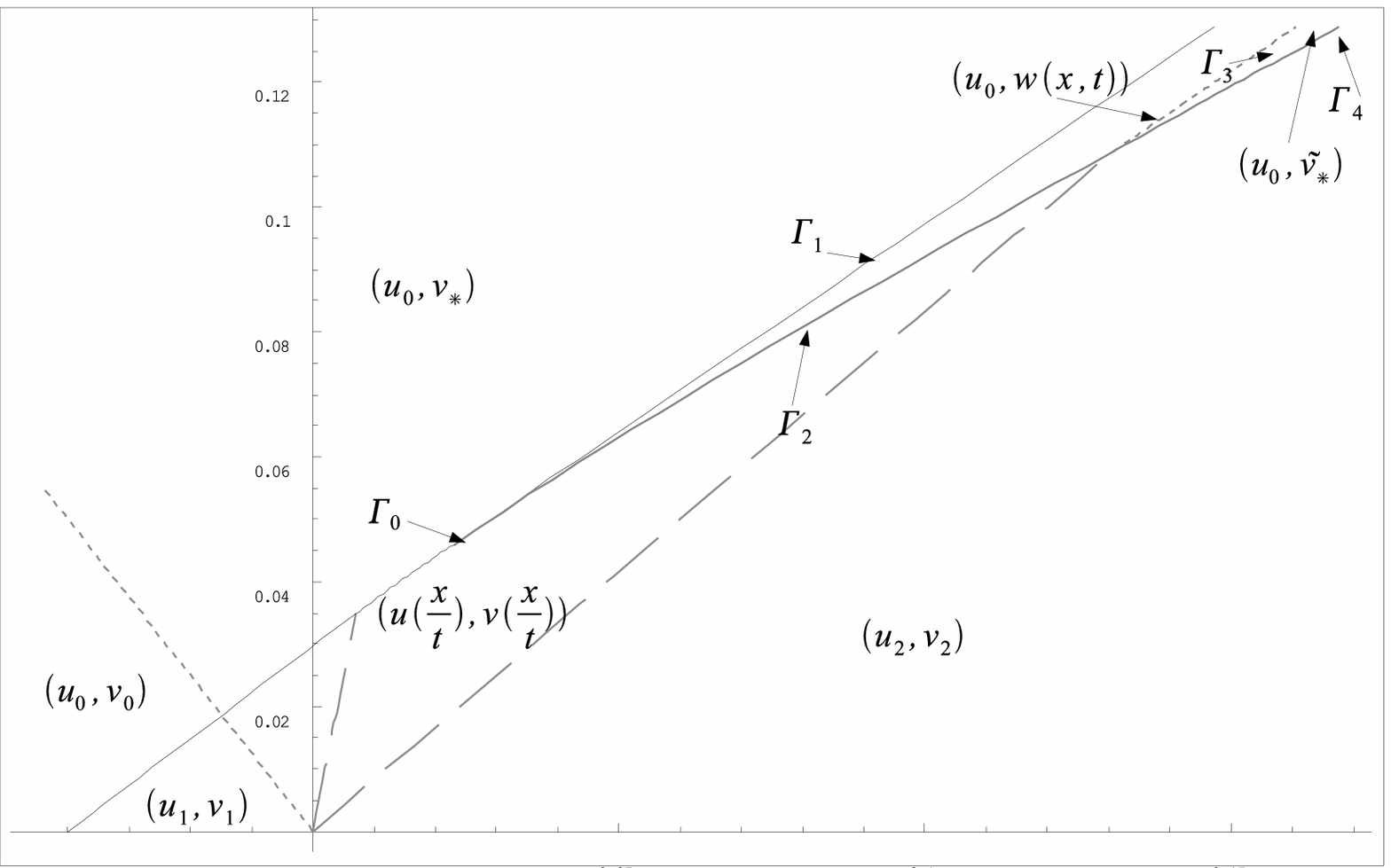}

{\small Fig.~5.}
\end{center}

\noindent
{\bf Case 5.}  Suppose that a delta shock wave starts from the point
$(0,a^{2})$, $a>0$ and meets a coupled pair of contact discontinuity
and rarefaction wave at some point $(x_{0},t_{0})$. This is possible
if $u_{0}<u_{1}$, $u_{1}>u_{2}+2$. Suppose that the rarefaction wave 
is centered (starts from $(0,0)$).

The point $(x_{0},t_{0})$ can be easily found by solving the equations
\begin{equation*}
\begin{split}
x-a^{2}={u_{1}+u_{2} \over 2}t, \;\; x=u_{1}t, \text{ i.e.} \\
t_{0}={2a^{2} \over u_{1}-u_{2}}, \;\; x_{0}={2a^{2}u_{1} \over u_{1}-u_{2}}.
\end{split}
\end{equation*}

In the beginning of the interaction of the rarefaction and the delta shock wave
the situation is quite similar to the one in the previous case. The solution 
is given by a delta shock wave supported by 
$\Gamma_{0}=\{ (t,c(t)):\;\; t>t_{0}\}$, where $c(t)$ is a solution to 
\begin{equation} \label{e9}
-c'(t)\Big( {c(t) \over t} - u_{2}\Big)
+{1\over 2} \Big( \Big({c(t)\over t}\Big)^{2}-u_{2}^{2}\Big)=0,\;\;
c(t_{0})=x_{0},
\end{equation}
i.e. 
\begin{equation*}
c(t)=u_{2}t+a\sqrt{2(u_{1}-u_{2})t},\;\; t>t_{0}.
\end{equation*}
Equation (\ref{e9}) is in fact Rankine-Hugoniot condition for (\ref{e1}).

The left- and right-hand side coefficients
of the new delta shock wave, $\alpha_{0}(t)$ 
and $\alpha_{1}(t)$, can be found in the same way as in the previous case.
If $(u(t),v(t))$ is the value of the rarefaction wave, then on the left-hand 
side of $\Gamma_{0}$ the new delta shock wave 
takes value $(u(t),v(t))|_{\Gamma_{0}}
=( c(t)/t,v_{1}\exp(c(t)/t-u_{1}))$ and on the right-hand side it equals
$(u_{2},v_{2})$. Only the overcompressibility condition is still in question.
The first condition for overcompressibility on the right-hand side
is always satisfied, since
$c'(t)=u_{2}+a\sqrt{2(u_{1}-u_{2})}/\sqrt{4t}>u_{2}$. For the 
overcompressibility it is necessary that also characteristic lines
run into the shock from the left-hand side.
\begin{equation*}
c'(t)=u_{2}+{a\sqrt{2(u_{1}-u_{2})} \over 2\sqrt{t}}
\leq u(c(t),t)-1=u_{2}+{a\sqrt{2(u_{1}-u_{2})} \over \sqrt{t}}-1,
\end{equation*}
i.e.\
\begin{equation*}
{a\sqrt{2(u_{1}-u_{2})} \over 2\sqrt{t}}\geq 1.
\end{equation*}
This is true until the time $t=t_{s}$, where 
\begin{equation*}
t_{s}={a^{2}(u_{1}-u_{2})\over 2}, \;\;
x_{s}=c(t_{s})=a^{2}(u_{1}-u_{2})\Big( {u_{2}\over 2}+1\Big).
\end{equation*}
Thus, $u_{s}=x_{s}/t_{s}=u_{2}+2$.

The first case: $u_{2}+2>u_{0}$.

Then the termination of overcompressibility takes place within the rarefaction
fan and again we are in a position to use the intuition 
behind Lemma \ref{l1}, i.e.\
to look for a solution consisting of a delta contact discontinuity
supported by a curve $\Gamma_{1}$ and a shock wave supported by some 
other curve $\Gamma_{2}$. $\Gamma_{1}$ should be below $\Gamma_{2}$.

$\Gamma_{1}$ is the characteristic line of the equation
\begin{equation*}
v_{t}+(u-1)v_{x}=0
\end{equation*}
passing trough $(x_{s},t_{s})$. Using the fact that $u=c_{1}(t)/t$ on
$\Gamma_{1}=\{(t,c_{1}(t)),\;\; t>t_{s}\}$, one can find such a function
$c_{1}$ by solving the initial value problem
\begin{equation*}
c_{1}'(t)={c_{1}(t)\over t} -1,\;\; c_{1}(t_{s})=x_{s}.
\end{equation*}
The unique solution to the above problem can be easily found
\begin{equation*}
c_{1}(t)=t\Big(-\log t +\log \Big({a^{2}(u_{1}-u_{2})\over 2}\Big)+u_{2}+2\Big).
\end{equation*}
Using (\ref{e1}) and the Rankine-Hugoniot condition, the curve
$\Gamma_{2}=\{(c_{2}(t),t)\; t>t_{s}\}$ is uniquely determined
by a solution to
\begin{equation*}
-c_{2}'(t)\Big(u_{2}-{c_{2}(t)\over t}\Big) 
+{1\over 2}\Big(u_{2}^{2}-\Big({c_{2}(t)\over t}\Big)^{2}\Big)=0,
\;\; c_{2}(t_{s})=x_{s},
\end{equation*}
i.e.\
\begin{equation*}
c_{2}(t)=u_{2}t+a\sqrt{2(u_{1}-u_{2})t}.
\end{equation*} 
One can see that it equals to the function $c(t)$ from the previous case.

One has to prove that $\Gamma_{1}$ is actually strictly below the curve 
$\Gamma_{2}$.

Since $c_{1}(t_{s})=c_{2}(t_{s})$ and $c_{1}'(t_{s})=c_{2}'(t_{s})=u_{s}-1$,
it is enough to compare $c_{1}'(t)$ and $c_{2}'(t)$, for $t>t_{s}$.
\begin{equation*}
\begin{split}
c_{1}'(t) & = u_{2}+1-\log t +\log \Big({a^{2}(u_{1}-u_{2})\over 2}\Big) \\
& = u_{2}+1+\log \Big({a^{2}(u_{1}-u_{2})\over 2t}\Big) \\
c_{2}'(t) & = u_{2}+{a\sqrt{u_{1}-u_{2}} \over \sqrt{2t}}.
\end{split}
\end{equation*}
$c_{2}'(t)> c_{1}'(t)$ if 
\begin{equation}\label{e10}
{a\sqrt{u_{1}-u_{2}} \over \sqrt{2t}}-1
> \log \Big({a^{2}(u_{1}-u_{2})\over 2t}\Big), \;\; t>t_{s}.
\end{equation}
But the last relation is true; one can check it by changing the variables,
and noticing that $\sqrt{y}-1 > \log (y)$, for $y\in (0,1)$.
\medskip

Denote by $A$ the region between 
$\Gamma_{1}$ and $\Gamma_{2}$ for $t>t_{s}$.

The value of $u$ is $u(x,t)=x/t$ inside $A$. Therefore, (\ref{e1}) 
is satisfied. Now, $\Gamma_{1}$ is the support of the delta contact 
discontinuity, and we are trying to find the value of $v$ in this area.

First, let $v_{\ast}(t)$ denote the value of $v$ on the left-hand side of
$\Gamma_{2}$. The value of $u$ there is given by
\begin{equation*}
u|_{\Gamma_{2}}={x\over t}|_{\Gamma_{2}}={c_{2}(t)\over t}=u_{2}
+{a\sqrt{2(u_{1}-u_{2})} \over \sqrt{t}}.
\end{equation*}
The values of $u$ and $v$ on the right-hand side of $\Gamma_{2}$ 
are $u_{2}$ and $v_{2}$, respectively. The Rankine-Hugoniot condition
for (\ref{e2}) gives
\begin{equation*}
-c_{2}'(t)(v_{2}-v_{\ast}(t))+\Big((u_{2}-1)v_{2}
-\Big({c_{2}(t)\over t}-1\Big)v_{\ast}(t)\Big)=0.
\end{equation*}
Solving the above equation, one gets
\begin{equation*}
v_{\ast}(t)={\sqrt{t}+B \over \sqrt{t}-B}v_{2},\;\; 
B={a\sqrt{2(u_{1}-u_{2})}\over 2}=\sqrt{t_{s}}.
\end{equation*}

Denote by $w(x,t)$ the value of $v$ inside $D$. Then $w$ is the solution to 
the linear partial differential equation
\begin{equation*}
w_{t}+\Big( {x\over t}-1\Big)w_{x}=0,\;\; w|_{\Gamma_{2}}=v_{\ast}(t).
\end{equation*}
The solution of the above equation is a constant along the characteristic
curves
\begin{equation*}
\gamma:={dx \over dt}={x\over t}-1, \text{ where }
\gamma|_{\Gamma_{2}} \text{ is known.}
\end{equation*}
In particular, $w$ tends to infinity near $\Gamma_{1}$, but in
locally integrable fashion because $v_{\ast}(t)={\mathcal O}(1/(\sqrt{t}
-\sqrt{t_{s}}))$ as $t \rightarrow t_{s}$.

Now, we shall look for an exit of the delta contact discontinuity and 
the shock wave supported by $\Gamma_{2}$ trough the rarefaction wave.

The line $x=u_{0}t$ and the curve $x=c_{1}(t)$ always has an interaction 
point, say $(\tilde{t},\tilde{x})$, for $u_{0}<u_{2}+2$. The equation
\begin{equation*}
u_{0}t=(u_{2}+2)t+t\log (B^{2}/t)=\log (t_{s}/t).
\end{equation*}
This equation has a unique solution $\tilde{t}>t_{s}$.

The next question is whether the curve $x=c_{2}(t)$ intersects the
line $x=u_{0}t$ or not. An intersection takes place, if the equation 
\begin{equation*}
u_{0}t=u_{2}t+2B\sqrt{t}, \text{i.e. }
u_{0}-u_{2}=2B/\sqrt{t}
\end{equation*}
has a solution $t>t_{s}$. If $u_{0}<u_{2}$, there is no solution. If 
$u_{0}>u_{2}$, then the solution $t=4B^{2}/(u_{0}-u_{2})^{2}=
4t_{s}/(u_{0}-u_{2})^{2}$ is bigger than $t_{s}$, because $2>u_{0}-u_{2}$.

In both cases we have to solve initial data problem 
for (\ref{e1},\ref{e2}), given by
\begin{equation*}
u=\begin{cases}u_{0}, & x<\tilde{x} \\ u_{0}, & x>\tilde{x}, \end{cases}
\;\; v=\left\{ \begin{aligned}v_{\ast}, & \;\; x<\tilde{x} \\ w(x,x/u_{0}), 
& \;\; x>\tilde{x} 
\end{aligned}\right\} + \gamma_{s}\delta_{(\tilde{t},\tilde{x})},
\end{equation*}
where $w$ is the right-hand side of $v$ in the region $A$, 
constant along the characteristics of
\begin{equation*}
w_{t}+\Big({x \over t} -1\Big) w_{x}=0, \;\; \gamma: {dx\over dt}
={x\over t} -1.
\end{equation*}
long the line $x=u_{0}t$ the slope of these characteristic curves is 
$u_{0}-1$. Thus we may continue the solution to the left of 
$x=u_{0}t$ as a delta contact discontinuity
\begin{equation*}
u(t,x)=u_{0},
\;\; v(t,x)=\left\{ \begin{aligned} 
v_{\ast}, & \;\; x-\tilde{x}<(u_{0}-1)(t-\tilde{t)} 
\\ w(t,x), & \;\; x-\tilde{x}>(u_{0}-1)(t-\tilde{t}) 
\end{aligned}\right\} + \gamma_{s}
\delta_{x-\tilde{x}=(u_{0}-1)(t-\tilde{t})},
\end{equation*}
where $w(t,x)$ is again constant along the lines with slope $u_{0}-1$. 
Denote by $\Gamma_{3}$ the line $x-\tilde{x}=(u_{0}-1)(t-\tilde{t})$.
There is no further intersection with the original contact discontinuity
along the (parallel) line $x=(u_{0}-1)t$.
In the case $u_{0}<u_{2}$, the solution is complete (See Fig.\ 6).

\begin{center}
\includegraphics*[scale=0.65]{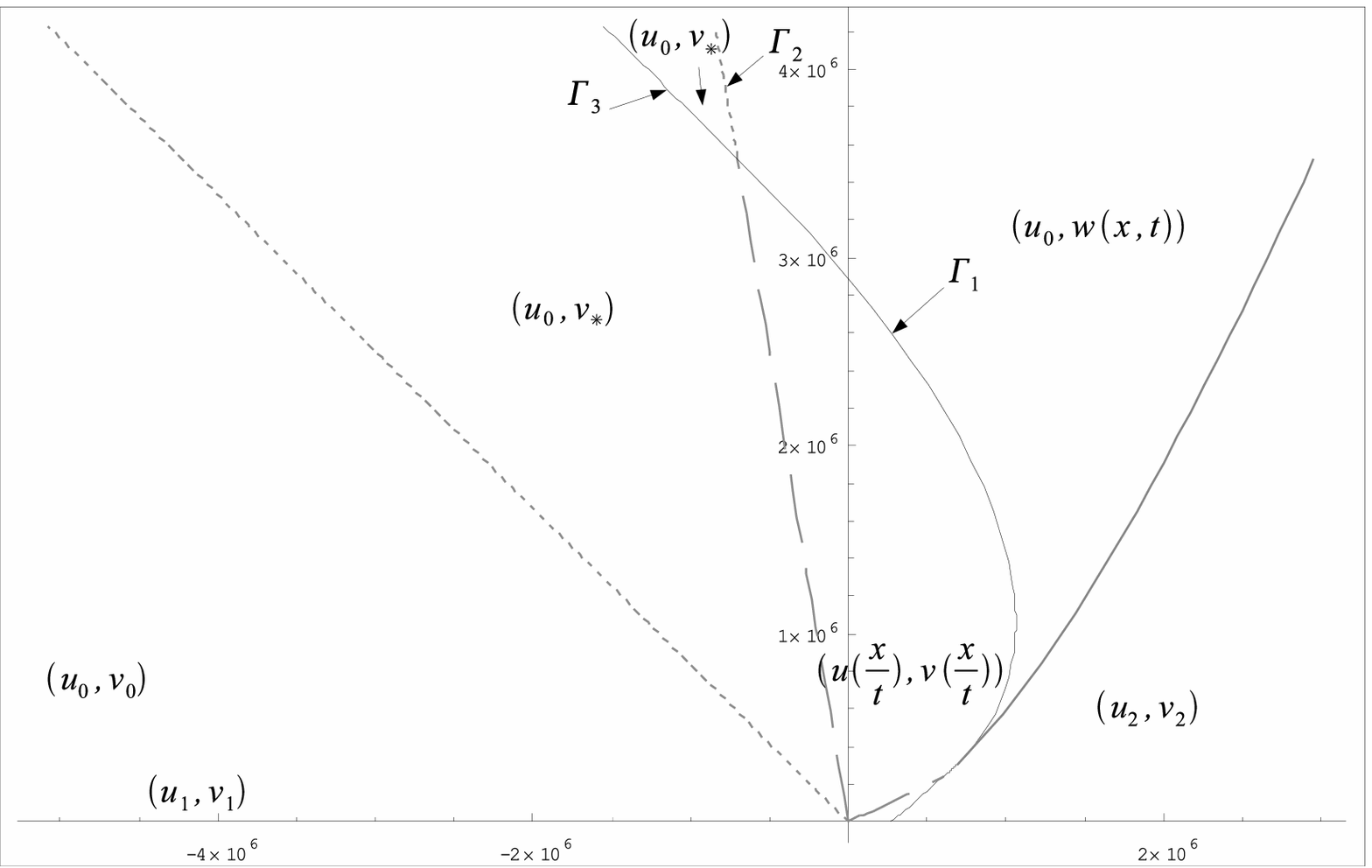}

{\small Fig.\ 6.}
\end{center}
\medskip

In the case $u_{0}>u_{2}$ we still have to consider the region above the 
intersection point $(\tilde{\tilde{x}},\tilde{\tilde{t}})$ of 
$x=c_{2}(t)$ with $x=u_{0}t$. In this case 
\begin{equation*}
u_{2}<u_{0}<u_{2}+2
\end{equation*}
and we can connect a constant left-hand state $u_{0}$ to the constant 
right-hand state $u_{2}$ by by a shock wave in $u$. 

This shock wave supported by $\Gamma_{5}$ 
has speed $(u_{0}+u_{2})/2$ and actually is tangent to the line
$x=c_{2}(t)$ at the intersection point $(\tilde{\tilde{x}},\tilde{\tilde{t}})$.
It follow a classical contact discontinuity starting from the point 
$(\tilde{\tilde{x}},\tilde{\tilde{t}})$ with speed $u_{0}-1$, 
supported by the line $\Gamma_{3}$ and this connects 
the region when $v=w(\xi,t)$ has been determined by the initial data 
along the line $x=u_{0}t$.  The value of $v$ between $\Gamma_{4}$
and $\Gamma_{5}$ is $\tilde{v}_{\ast}=(2+u_{0}-u_{2})/(2+u_{2}-u_{0})$.
(See Fig.\ 7).

\begin{center}
\includegraphics*[scale=0.65]{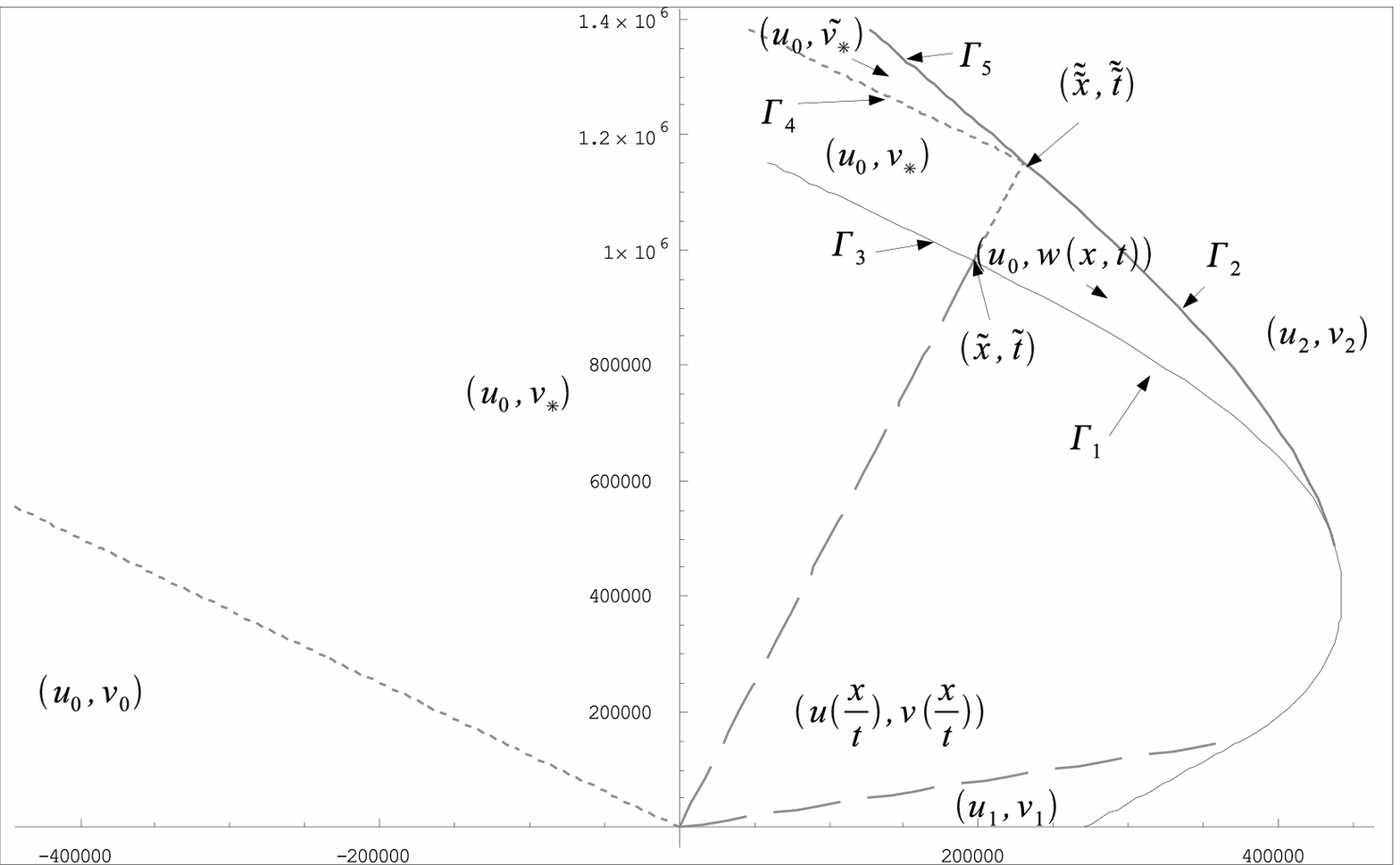}

{\small Fig.~7.}
\end{center}
\medskip
\medskip
\medskip

The second case is $u_{0}>u_{s}=u_{2}+2$. Then there is no bifurcation of 
the delta contact discontinuity supported by $\Gamma_{0}$. After 
$\Gamma_{0}$ intersects the line $x=u_{0}t$ at $(t_{1},x_{1})$, say,
the solution can be continued into the region 
$(u_{0}-1)t<x<u_{0}t$ by a simple delta contact discontinuity on the 
line $\Gamma_{3}:x-x_{1}=(u_{0}-1)(t-t_{1})$, where $u$ has the 
constant value $u_{0}$ and $v$ has the value $v_{\ast}$ and $v_{2}$
on the left and right-hand side, respectively, and a constant strength 
delta function is placed on the line $\Gamma_{3}$.
This concludes investigation of all possible cases.


\begin{thebibliography}{99}

\bibitem{Bre} Bressan, A., Hyperbolic Systems of Conservation Laws,
Preprins S.I.S.S.A., Trieste, Italy. 

\bibitem{c}  Colombeau, J.\ F., Elementary Introduction in New
Generalized Functions North Holland, 1985.

\bibitem{hlf} Hayes, B.\ T.\ and  Le Floch, P.\ G., 
'Measure solutions to a strictly hyperbolic system of
conservation laws',  {\it Nonlinearity} {\bf 9}, 1547-1563 (1996).

\bibitem{k} Keyfitz, B.\ L., 'Conservation laws, delta shocks and
singular shocks', In: M.\ Grosser  et al: Nonlinear Theory of
Generalized Functions, Research Notes in Math., Champman Hall/CRC, 1999.

\bibitem{kk}  Keyfitz, B.\ L.\ and  Kranzer, H.\ C.,
'Spaces of weighted measures for conservation laws with singular
shock solutions', {\it J.\ Diff.\ Eq.} {\bf 118},2, 420-451 (1995).

\bibitem{kor}  Korchinski, D.\ J., Solution of a Riemann Problem for
a $2\times 2$ System of Conservation Laws Possessing No Classical
Weak Solution, PhD Thesis, Adelphi University, Garden City, New York, 
1977.

\bibitem{lax}  P.\ D.\ Lax, Hyperbolic Systems of 
Conservation Laws and the Mathematical Theory of
Shock Waves, SIAM, Philadelphia, 1973.

\bibitem{mn1} Nedeljkov, M., 'Delta and singular delta locus for 
one dimensional systems of conservation laws', {\it Math.\ Meth.\ Appl.\ Sci.}
{\bf 27}, 931--955 (2004).

\bibitem{mn2} Nedeljkov, M., 'Second delta locus - interaction of 
delta shock with shock waves', {\it Preprint}.

\bibitem{o}  Oberguggenberger, M., Multiplication of Distributions
and Applications to Partial Differential Equations,
Pitman Res.\ Not.\ Math.\ 259, Longman Sci.\ Techn.,  Essex, 
1992.

\bibitem{ow} Oberguggenberger, M.\ and  Wang, Y-G., 
'Generalized solutions to conservation laws',
{\it Zeitschr.\ Anal.\ Anw.} {\bf 13}, 7-18 (1994).

\bibitem{tzz}  Tan, D., Zhang, T.\ and  Zheng, Y., 
'Delta-shock waves as limits of vanishing viscosity for 
hyperbolic systems of conservation laws',
{\it J.\ Diff.\ Eq.} {\bf 112}, 1-32 (1994).

\end{thebibliography}
\end{document}